\documentclass[10pt,twoside]{article}

\usepackage{amssymb}
\usepackage{amsmath}
\usepackage{amsthm}

\usepackage{mathrsfs}
\usepackage[T1]{fontenc}            %% Polish fonts
\usepackage[utf8x]{inputenc}        %% File encoding
\usepackage[polish,british]{babel}  %% Languages (for hyphenation rules), british is prefered

\usepackage{a4}                     %% paper layout for A4-paper.
\usepackage{verbatim}               %% I use verbatim mainly for comment enviroment.
\usepackage[all,line]{xy}                %% xy is for diagrams
\usepackage[dvips]{graphicx}        %% graphicx is for inserting eps or other figures in the text.
\usepackage{url}                    %% for producing links, mainly in the references.

\usepackage[textsize=tiny]{todonotes}
\def\CC{{\mathbb C}}

\def\PP{{\mathbb P}}

\def\ZZ{{\mathbb Z}}

\newcommand{\ccE}{{\mathcal{E}}}

\newcommand{\ccH}{{\mathcal{H}}}

\newcommand{\ccL}{{\mathcal{L}}}
\newcommand{\ccM}{{\mathcal{M}}}

\newcommand{\ccO}{{\mathcal{O}}}
\newcommand{\ccP}{{\mathcal{P}}}

\newcommand{\ccS}{{\mathcal{S}}}

\newcommand{\ccU}{{\mathcal{U}}}

\def\gotm{\mathfrak{m}}
\newcommand{\gote}{\mathfrak{e}}
\newcommand{\gotf}{\mathfrak{f}}
\newcommand{\gotg}{\mathfrak{g}}
\newcommand{\gotp}{\mathfrak{p}}
\newcommand{\gotso}{\mathfrak{so}}

\def\half{\frac{1}{2}}

\def\ud{\mathrm{d}}

\DeclareMathOperator{\codim}{codim}

\DeclareMathOperator{\Hilb}{Hilb}
\DeclareMathOperator{\Hom}{Hom}

\DeclareMathOperator{\id}{id}

\DeclareMathOperator{\Pic}{Pic}
\DeclareMathOperator{\Proj}{Proj}

\DeclareMathOperator{\sspan}{span}

\newcommand{\SheafyProj}{\ccP{}\mathit{roj}}

\def\Pic{\operatorname{Pic}}

\def\codim{\operatorname{codim}}

\newcommand{\reduced}[1]{{#1}_{\operatorname{red}}}

\newcommand{\set}[1]{\left\{#1\right\}}

\def\Wedge#1{{\textstyle{\bigwedge\nolimits}^{\! #1}}}

\numberwithin{equation}{section}
\newtheorem{thm}[equation]{Theorem}
\newtheorem{prop}[equation]{Proposition}
\newtheorem{lemma}[equation]{Lemma}
\newtheorem{cor}[equation]{Corollary}
\newtheorem{conjecture}[equation]{Conjecture}
\theoremstyle{remark}
\newtheorem{example}[equation]{Example}
\newtheorem{rem}[equation]{Remark}

\theoremstyle{definition}
\newtheorem{defin}[equation]{Definition}

\newtheorem{notation}[equation]{Notation}
\newtheorem{setting}[equation]{Setting}

\newenvironment{prf}[1][]
  {\medskip\par\noindent{\bf Proof#1. }}
  {\nopagebreak\par\rightline{$\Box$} \medskip}

\def\noprf{\nopagebreak\par\rightline{$\Box$} \medskip}

\newcommand{\Buczynski}{Buczy\'n{}ski}
\newcommand{\Jaroslaw}{Jaros\l{}aw}
\newcommand{\JaBu}{\Jaroslaw{} \Buczynski{}}
\newcommand{\Wisniewski}{Wi\'s{}\-nie\-wski}

\author{\JaBu}

\title{Duality and integrability on contact Fano manifolds}

 \date{15th February 2021}
\newcommand{\ContactDistributionSequence}[1]{%
\[
0 \to F \to T #1 \stackrel{\theta}{\to} L \to 0
\]
}

\newcommand{\divisor}[1]{{#1}^{div}}
\newcommand{\point}[1]{{#1}^{\PP}}
\newcommand{\hyperplane}[1]{{#1}^{\PP \perp}}

\newcommand{\LinSyst}{{\left\langle D \right\rangle}}
\newcommand{\LinSystDual}{{\LinSyst^*}}

\newcommand{\component}[1]{{(#1)^{\bullet}}}
\newcommand{\legPF}{\mathscr{C}}

\begin{document}
\maketitle

\begin{abstract} 
We address the problem of classification of contact Fano manifolds.
It is conjectured that every such manifold is necessarily homogeneous.
We prove that the Killing form,  
the Lie algebra grading and parts of the Lie bracket can be read from geometry of an arbitrary
contact manifold. 
Minimal rational curves on contact manifolds (or contact lines) 
and their chains are the essential ingredients for our constructions.
Along the way we collect several facts about $\PP^1$-bundles admitting a contractible section, which might be of independent interest.
\end{abstract}

\medskip
{\footnotesize
\noindent\textbf{author's e-mail:} jabu@mimuw.edu.pl

\noindent\textbf{keywords:} 
complex contact manifold, Fano variety, minimal rational curves,
adjoint variety, Killing form, Lie bracket, Lie algebra grading, 
complex contact manifold;

\noindent\textbf{AMS Mathematical Subject Classification 2010:}
Primary: 14M17; Secondary: 53C26, 14M20, 14J45;

\medskip

This version of the article consolidates the published version \cite{jabu_contact_duality_and_integrability} 
and its erratum \cite{jabu_contact_duality_and_integrability_errata}.}

\section{Introduction}\label{section-introduction}

In this article we are interested in the classification of contact Fano manifolds.
We review the relevant definitions in \S\ref{section-preliminaries}.
So far the only known examples of contact Fano manifolds are obtained as follows.
For a simple Lie group $G$ consider its adjoint action on $\PP(\gotg)$,
where $\gotg$ is the Lie algebra of $G$. 
This action has a unique closed orbit $X$ and this $X$ has a natural contact structure.
In this situation $X$ is called a \emph{projectivised minimal nilpotent orbit},
or the \emph{adjoint variety} of $G$.
By the duality determined by the Killing form, 
equivalently we can consider the coadjoint action of $G$ on $\PP(\gotg^*)$ and $X$ is isomorphic 
to the unique closed orbit in $\PP(\gotg^*)$.

\begin{conjecture}[LeBrun, Salamon] \label{con_contact}
If $X$ is a Fano complex contact manifold, then $X$ is 
the adjoint variety of a simple Lie group $G$.
\end{conjecture}

This problem originated with a problem in Riemannian geometry.
In \cite{berger}  a list of all possible holonomy groups 
of simply connected Riemannian manifolds is given. 
The existence problem for all the cases has been solved locally. 
Also compact simply connected non-homogeneous examples with most of the possible holonomy groups were constructed
with the unique exception of the quaternion-K\"{a}hler manifolds.
It is conjectured that the compact simply connected quaternion-K\"{a}hler manifolds must be symmetric spaces, especially those with positive scalar curvature
(see \cite{lebrun} and references therein).

\begin{conjecture}[LeBrun, Salamon]\label{con_qK}
Let $M$ be a compact simply-connected qua\-ter\-nion-K\"{a}hler manifold with positive scalar curvature.
Then $M$ is a homogeneous symmetric space 
(more precisely, it is one of the Wolf spaces  --- see \cite{wolf}).
\end{conjecture}

The relation between the conjectures is given by the construction of a twistor space.
The $S^2$-bundle of complex structures on tangent spaces to a quaternion-K\"{a}hler manifold $M$ 
 is called the twistor space of $M$.
 If $M$ is compact and has positive scalar curvature,
 and then the twistor space $X$ has a natural complex structure and 
 is a contact Fano manifold with a K\"{a}hler-Einstein metric. 
 In particular, the twistor space of a Wolf space is an adjoint variety.
Hence Conjecture~\ref{con_contact} implies Conjecture~\ref{con_qK}.
Conversely, if $X$ is a contact Fano manifold with 
 K\"{a}hler-Einstein metric, then 
 it is a twistor space of a quaternion-K\"{a}hler manifold with positive scalar curvature --- see \cite{lebrun}. 
 See also the \cite[\S4]{jabu_moreno_contact_survey} for a more detailed survey and references.

In order to study the non-homogeneous contact Fano manifolds (potentially non-existent) it is natural to assume
$\Pic X  \simeq \ZZ$ and further that $X$ is not isomorphic to a projective space. 
This only exludes the adjoint varieties of types $A$ and $C$
(see \S\ref{section-preliminaries} for more details).

With this assumption, 
we take a closer look at three pieces of the homogeneous structure on adjoint varieties:
the Killing form $B$ on $\gotg$, the Lie algebra grading
$\gotg = \gotg_{-2} \oplus \gotg_{-1} \oplus \gotg_{0} \oplus \gotg_{1} \oplus \gotg_{2}$
(see \cite[\S6.1]{landsbergmanivel02} and references therein)
and a part of the Lie bracket on $\gotg$.
Understanding the underlying geometry allows us to define the appropriate generalisations of these notions 
on arbitrary contact Fano manifolds.

An essential building block for our constructions is the notion of a \emph{contact line} 
(or simply \emph{line})
on $X$.
These contact lines were studied by Kebekus \cite{kebekus_lines1}, \cite{kebekus_lines2} and
\Wisniewski{} \cite{wisniewski}.
Also they are an instance of minimal rational curves, which are studied extensively.
The geometry of contact lines was the original motivation to study Legendrian subvarieties in projective space
(see \cite{jabu_dr} for an overview and many details).
We briefly review the subject of lines on contact Fano manifolds in \S\ref{section-review-lines}.
 
The key ingedient is the construction of a family of divisors $D_x$ parametrised by points $x \in X$
(see \S\ref{section-define-Dx}).
These divisors are swept by pairs of intersecting contact lines, one of which passes through $x$.
In other words, set theoretically $D_x$ is the set of points of $X$,
which can be joined with $x$ using at most $2$ intersecting contact lines.
The idea to study these loci comes from \Wisniewski{} \cite{wisniewski}
where he observed, that (under an additional minor assumption)
 these loci contain some non-trivial divisorial components
and he studied the intersection numbers of certain curves on $X$ with the divisorial components.
Here we prove all the components of $D_x$ are divisorial and draw conclusions from that observation
going into a different direction than those of \cite{wisniewski}.

\begin{thm}\label{thm-exist-Dx-and-duality}
   Let $X$ be a contact Fano manifold with $\Pic X \simeq \ZZ $ 
   and assume $X$ is not isomorphic to a projective space.
   Then the locus $D_x \subset X$ 
   swept by the pairs of intersecting contact lines, one of which passes through $x \in X$
   is of pure codimension $1$ and thus $D_x$ determines a divisor on $X$.
   Let $\LinSyst \subset H^0(\ccO(D_x))$ be the linear system spanned by these divisors.
   Let $\phi\colon X \to \PP\LinSystDual$ be the map determined by the linear system $\LinSyst$
   and let $\psi\colon X \to \PP\LinSyst$ be the map $x\mapsto D_x$.
   Then:
   \begin{enumerate}
     \item both $\phi$ and $\psi$ are regular maps.
     \item there exists a unique up to scalar non-degenerate bilinear form $B$ on $\LinSyst$,
           which determines an isomorphism $\PP\LinSystDual \simeq \PP \LinSyst$ 
	   making the following diagram commutative:
	   \[
              \xymatrix@C1.5pc@R0.5pc{
              & & & \PP \LinSystDual \ar[dd]^{\simeq} \\
              X\ar[urrr]^{\phi} \ar[drrr]^{\psi} && \\
              & & & \PP \LinSyst.}
	   \]
     \item The bilinear form $B$ is either symmetric or skew-symmetric.
     \item If $X \subset \PP(\gotg^*)$ is the adjoint variety of simple Lie group $G$,
          then $\LinSyst = \gotg$ and $B$ is the Killing form on $\gotg$.
   \end{enumerate}
\end{thm}

With the notation of the theorem,
after fixing a pair of general points $x, w \in X$ 
there are certain natural linear subspaces of $\LinSyst$,
which we denote $\LinSyst_{-2}$, $\LinSyst_{-1}$, $\LinSyst_{0}$, $\LinSyst_{1}$ and $\LinSyst_{2}$
(see \S\ref{section-grading} for details).
\begin{thm} \label{thm-grading-for-homogeneous}
  If $X \subset \PP(\gotg^*)$ is the adjoint variety of a simple Lie group $G$ with $\Pic X \simeq Z$ 
  and $X$ not isomorphic to a projective space,
  then there exists a choice of a maximal torus of $G$
  and a choice of order of roots of $\gotg$,
  such that 
  $\LinSyst_{i} = \gotg_{i}$ for every $i \in \set{-2,-1,0,1,2}$,
  where $\gotg = \gotg_{-2} \oplus \gotg_{-1} \oplus \gotg_{0} \oplus \gotg_{1} \oplus \gotg_{2}$
  is the Lie algebra grading of $\gotg$.
\end{thm}
 
Finally, if $X$ is the adjoint variety of $G$,
then there is a rational map 
\[
  [\cdot, \cdot]\colon X \times X \dashrightarrow \PP(\gotg),
\]
which is the Lie bracket on $\gotg$ (up to projectivisation).
Also there is a divisor $D \subset X \times X$,
such that for general $(x,z) \in D$ 
the Lie bracket $[x,z]$ is in $X$.
We recover this bracket restricted to $D$  for general contact manifolds:
\begin{thm} \label{thm-bracket}
  For $X$ and $D_x$ as in Theorem \ref{thm-exist-Dx-and-duality},
  let $D \subset X \times X$ be the divisor consisting of pairs $(x,z) \in X \times X$,
  such that $z \in D_x$.
  There exists a rational map $[\cdot , \cdot ]^D\colon D \dashrightarrow X$,
  such that $[x,z]^D = y$,
  where $y$ is an intersection point of a pair of contact lines that join $x$ and $z$.
  In particular, this intersection point $y$ and the pair of lines are unique for general pair $(x,z)$ $ \in D$.
  Moreover, if $X$ is the adjoint variety of a simple Lie group $G$, then $[ \cdot, \cdot]^D$ is the restricion of the Lie bracket.
\end{thm}

In \S\ref{section-preliminaries} we introduce and motivate our assumptions and notation.
In \S\ref{sec_bundles_of_P1s} we discuss $\PP^1$-bundles with a contractible section. This subsection comes from the erratum \cite{jabu_contact_duality_and_integrability_errata} 
and is required to clarify results about contact lines from \cite{kebekus_lines2}, see Theorem~\ref{thm-properties-of-Cx}.
It can also serve of independent interest.

In \S\ref{section-loci-swept-by-lines} we review the notion of contact lines and their properties.
We continue by studying certain types of loci swept by those lines and calculate their dimensions.
In particular we prove there Theorem~\ref{thm-Dx-is-a-divisor},
which is a part of results summarised in Theorem~\ref{thm-exist-Dx-and-duality}.
We also study the tangent bundle to $D_x$ as a subspace of $TX$.

In \S\ref{section-duality} we study the duality of maps $\phi$ and $\psi$
introduced in Theorem~\ref{thm-exist-Dx-and-duality} together with the consequences of this duality.
This section is culminated with the proof of Theorem~\ref{thm-exist-Dx-and-duality}.
 
In \S\ref{section-grading} we generalise the Lie algebra grading to arbitrary contact manifolds 
and prove Theorem~\ref{thm-grading-for-homogeneous}.

In \S\ref{section-cointegrable} we prove that certain lines are integrable with respect to 
a special distribution on $D_x$ and we apply this to prove Theorem~\ref{thm-bracket}.

\subsection*{Acknowledgements}

The author was supported by Marie Curie International Outgoing Fellowship. While preparing the corrigendum, the author was supported by the National Scienc Center of Poland (NCN) project ``Complex contact manifolds and geometry of secants'', 2017/26/E/ST1/00231.

The author would like to thank the following people for many enlightening 
discussions:
Jun-Muk Hwang,
Stefan Kebekus,
Joseph M.~Landsberg,
\Jaroslaw{} \Wisniewski{}
and Fyodor L.~Zak.
Comments of Laurent Manivel helped to improve the presentation in the paper.
In addition, the author is grateful to Joachim Jelisiejew, Micha{\l} Kapustka for further suggestions concerning the erratum.

\section{Preliminaries}\label{section-preliminaries}

Throughout the paper all our projectivisations $\PP$ are naive.
This means, if $V$ is a vector space, then $\PP V = (V \setminus 0)/ \CC^*$,
and similarly for vector bundles.

A complex manifold $X$ of dimension $2n+1$ is \emph{contact}
if there exists a vector subbundle $F \subset TX$ of rank $2n$ fitting into an exact sequence: 
\ContactDistributionSequence{X}
such that the derivative $\ud \theta \in H^0(\Wedge{2} F^* \otimes L)$
of the twisted form $\theta \in H^0(T^*X \otimes L)$ is nowhere degenerate.
In particular, $\ud \theta_x$ is a symplectic form on the fibre of contact distribution $F_x$.
See \cite[\S{}E.3 and Chapter~C]{jabu_dr} and references therein for an overview of the subject.

A projective manifold $X$ is Fano, if its anticanonical divisor ${K_X}^*=\Wedge{\dim X} TX$ is ample.

If $X$ is a projective contact manifold,
then by Theorem of Kebekus, Peternell, Sommese and \Wisniewski{} \cite{4authors}
combined with a result by Demailly \cite{demailly},
$X$ is either a projectivisation of a cotangent bundle 
to a smooth projective manifold 
or $X$ is a contact Fano manifold, with $\Pic X \simeq \ZZ$.
In the second case, since $K_X \simeq (L^*)^{\otimes(n+1)}$, by \cite{kobayashi_ochiai},
either $X \simeq \PP^{2n+1}$ or $\Pic X = \ZZ \cdot [L]$.
Here we are interested in the case $X \not \simeq \PP^{2n+1}$.
Thus our assumption spelled out below only exclude some well understood cases 
(the projectivised cotangent bundles and the projective space)
and they agree with the assumptions of 
Theorems~\ref{thm-exist-Dx-and-duality}, \ref{thm-grading-for-homogeneous} and~\ref{thm-bracket}.
\begin{notation}\label{notation-with-assumptions-about-contact-mflds}
   Throughout the paper $X$ denotes a contact Fano manifold with $\Pic X$ generated by the class of $L$,
   where $L= TX/F$ and $F\subset TX$ is the contact distribution on $X$. 
   We also assume $\dim X = 2n +1$.
\end{notation}
From Theorem of Ye \cite{ye} it follows that $n \ge 2$.

We will also consider the homogeneous examples of contact manifolds (i.e. the adjoint varieties).
Thus we fix notation for the Lie group and its Lie algebra.

\begin{notation}\label{notation-with-assumptions-about-group}
   Throughout the paper $G$ denotes a simple complex Lie group, not of types $A$ or $C$ 
   (i.e.~not isomorphic to $SL_n$ nor $Sp_{2n}$ nor their discrete quotients).
   Further $\gotg$ is the Lie algebra of $G$.
   Thus $\gotg$ is one of $\gotso_n$ (types $B$ and $D$),
   or one of the exceptional Lie algebras 
   $\gotg_2$, $\gotf_4$, $\gote_6$, $\gote_7$ or $\gote_8$.
\end{notation}

The contact structure on $\PP^{2n-1} = \PP(\CC^{2n})$ is determined by a symplectic form 
$\omega$ on $\CC^{2n}$. 
 The precise relation between the contact and symplectic structures is decribed for instance in 
\cite[\S{}E.1]{jabu_dr} (see also \cite[Ex.~2.1]{lebrun}).
In particular, for all $x \in X$, the projectivisation of a fibre of the contact distribution $\PP F_x$
comes with a natural contact structure.

Let $M$ be a projective contact manifold
(in our case $M = X$ with $X$ as in Notation~\ref{notation-with-assumptions-about-contact-mflds}
 or $M = \PP^{2n-1}$).
A subvariety $Z \subset M$ is \emph{Legendrian},
 if for all smooth points $z \in Z$ the tangent space $T_z Z$
 is contained in the contact distribution of $M$ and $Z$ is of pure dimension $\half(\dim M -1)$.

Recall from \cite[Lecture~20]{harris} or \cite[III.\S3,\S4]{mumford} the notion of \emph{tangent cone}.
For a subvariety $Z \subset X$, and a point $x \in Z$ let $\tau_x Z \subset T_x X$ 
be the tangent cone of $Z$ at $x$.
In this article we will only need the following elementary properties of the tangent cone:
\begin{itemize}
 \item $\tau_x Z$ is an affine cone (i.e.~it is invariant under the standard action of $\CC^*$ on $T_x X$).
 \item $\dim_x Z = \dim \tau_x Z$ and thus if $Z$ is irreducible, then $\dim Z = \dim \tau_x Z$.
 \item If $x \in Z_1 \subset Z_2 \subset X$, then  $\tau_x Z_1 \subset \tau_x Z_2$.
 \item If $Z$ is smooth at $x$, then $\tau_x Z = T_x Z$. 
\end{itemize}
Since $\tau_x Z$ is a cone, let $\PP \tau_x Z \subset \PP T_x X$ be the corresponding projective variety.

\subsection{Bundles of projective lines}\label{sec_bundles_of_P1s}

The situation we are primarily interested in is the following: 
\begin{setting}\label{set_bundle_with_contraction}
Let $Y, U, B$ be normal projective varieties, such that 
$\pi \colon U \to B$ is an analytically locally trivial $\PP^1$-bundle with a section $s\colon B \to U$ and a map $\xi\colon U \to Y$ that is birational onto its image and the section $s(B)$ is contracted to a point $y\in Y$:
\[
   \xymatrix{
      U \ar[r]^{\xi}\ar[d]^{\pi} & Y\\
      B. \ar@/^/[u]^{s} 
   }
\]
\end{setting}

The lemmas in this subsection seem to be classical and known to experts, but hard to reference explicitly.
Setting~\ref{set_bundle_with_contraction} (with more details spelled out as Setting~\ref{set_reducible_base_immersive_into_blup}) appears in \cite[\S6.1]{kebekus_lines2}. It seems that most of the statements in this subsection have been known to Kebekus, and have been implicitly used in his arguments. 
However, since there is a gap in his arguments (more on this below), we provide all necessary details in order to avoid analogous bugs at the cost of being somewhat tedious. 

In the first two lemmas we focus on $\PP^1$-bundles with one section.

\begin{lemma}\label{lem_projectivisation_of_vb}
   Suppose $U$ and $B$ are normal projective varieties and $\pi \colon U \to B$ is an analytically locally trivial $\PP^1$-bundle with a section $s \colon B \to U$.
   Then $U = \PP(\ccE)$, where $\ccE$ is a vector bundle of rank $2$ over $B$. In particular, $U$ is a Zariski locally trivial $\PP^1$-bundle over $B$.
   Moreover, 
   \[
      \pi_*(k \cdot s(B)) = \operatorname{Sym}^k (\ccE \otimes \ccM^*),
   \]
   where $k\cdot s(B)$ is the divisor on $U$ 
   (the $k$-fold multiple of $s(B)$) and $\ccM \subset \ccE$ is the line subbundle corresponding to $s(B)$.
\end{lemma}
\begin{prf}
   The image $s(B)$ is a divisor on $U$ which is analytically locally given as vanishing of a single holomorphic function. Since $U$ is projective, by GAGA, $s(B)$ is a Cartier divisor. 
   Let $\ccO_{U}(s(B))$ be the corresponding line bundle.
   Note that $\ccO_U(s(B))$ is $\pi$-ample by \cite[Thm~1.7.8]{lazarsfeld}, and its restriction to any fibre $\PP^1$ is $\ccO_{\PP^1}(1)$. 
   Thus by \cite[Prop.~4.10]{araujo_druel_codim1_del_Pezzo_foliations} the bundle is a projectivisation of a vector bundle $\ccE$ as claimed. 
\end{prf}

Note that in the setting of 
   Lemma~\ref{lem_projectivisation_of_vb} the vector bundle $\ccE$ such that 
   $U \simeq \PP(\ccE)$ 
   is not unique. 
   Instead, $\ccE$ is well defined 
   up to a twist by a line bundle on $B$.
Nevertheless, the formula 
   $\pi_*(k \cdot s(B)) 
   = \operatorname{Sym}^k (\ccE \otimes \ccM^*)$
   is independent of the choice of $\ccE$.
   
\begin{lemma}\label{lem_scheme_structures_are_equal}
   Suppose $U$ and $B$ are reduced schemes and $\pi \colon U \to B$ is a projective morphism that is a Zariski locally trivial $\PP^1$-bundle with a section $s \colon B \to U$.
   Suppose further that $S\subset U$ is a subscheme $S$
     supported on $s(B)$.
   Then $S= s(B)$ (as schemes) if and only if for all $b\in B$ 
     the intersection of the fibre $\pi^{-1}(b)$
     and $S$ is the single reduced point $\set{s(B)}$.
\end{lemma}
\begin{prf}
   If $S=s(B)$, then the claim on intersection is clear. 
   So suppose that $S$ intersects a fibre $\pi^{-1}(b)$
     in the single reduced point. 
   We will show, that $S$ agrees with $s(B)$ in some neighbourhood of $s(b)$.
   
   Let $\ccO_{S,s(b)}$ be the local ring of $S$ near $s(b)$ and let 
   $I\subset \ccO_{S,s(b)}$ be the ideal of $s(B)$.
   Note that $\ccO_{S,s(b)}/I \simeq \ccO_{B,b}$
      and by local triviality also 
      $\ccO_{S,s(b)} =\ccO_{B,b} \oplus I$ as $\ccO_{B,b}$-modules.
   Let $\gotm_b \subset \ccO_{B,b}$ be the maximal ideal.
   The assumption on the intersection is translated into algebra as
      \[
         \ccO_{S,s(b)} \otimes_{\ccO_{B,b}} \ccO_{B,b}/\gotm_b 
         = \ccO_{B,b}/\gotm_b.
      \]
   Thus $ I \otimes_{\ccO_{B,b}} \ccO_{B,b}/\gotm_b =0$ 
      and by Nakayama's Lemma $I=0$ as claimed.
\end{prf}

In the following lemma we focus on $\PP^1$-bundles that admit two disjoint sections.   
   
\begin{lemma}\label{lem_P1_bundle_with_two_sections}
   Suppose $U$ and $B$ are normal projective varieties and $\pi \colon U \to B$ 
     is an analytically locally trivial $\PP^1$-bundle with a section $s \colon B \to U$.
   Suppose moreover that there exists another section $s'\colon B\to U$ whose image is disjoint from $s(B)$.
   Then
   \begin{enumerate}
    \item \label{item_Pic_of_total_space_of_bundle}
     $\Pic U = \pi^*\Pic(B) \oplus \ZZ \cdot s'(B)$ and the pullback $s^*\colon \Pic(U)\to \Pic(B)$ 
    is the projection on the first component $\pi^*\Pic(B) \simeq \Pic(B)$, and
    \item \label{item_total_space_as_decomposable_bundle}
    $U=\PP(\ccO_B \oplus \ccM)$, where $\ccM$ is a line bundle on $B$ such that the pullback 
      of the corresponding Cartier divisor $\pi^*\ccM$ is linearly equivalent to the divisor $s(B)-s'(B)$.
   \end{enumerate}
\end{lemma}
\begin{prf}
   $U = \PP(\ccE)$ by Lemma~\ref{lem_projectivisation_of_vb},
      and by the assumptions $\ccE$ admits two line subbundles which 
      do not intersect away from the zero section.
   Thus $\ccE$ is decomposable, and twisting by one of the components, 
     without loss of generality we can assume $\ccE = \ccO_B \oplus \ccM$ for some line bundle $\ccM$, where $\PP(\ccO_B) \subset \ccU$ is the image of the section $s$, while $\PP(\ccM)$
     is the image of $s'$.
     
   In this setting~\ref{item_Pic_of_total_space_of_bundle}
   is straightforward.
   To finish the 
   proof of~\ref{item_total_space_as_decomposable_bundle}
   we must show $\pi^*\ccM =  s(B)-s'(B)$ in $\Pic(U)$. 
   Since $(s(B)-s'(B)). \PP^1 = 0$ (where $\PP^1$ is the general fibre of $\pi$) 
   by~\ref{item_Pic_of_total_space_of_bundle}
   we must have $s(B)-s'(B) \in \pi^*\Pic (B)$, 
   thus:
   \[
     s^*(s(B)-s'(B)) = s^*(s(B)) = N_{s(B)\subset U} =
     N_{\PP(\ccO_B)\subset \PP(\ccO_B \oplus \ccM)} = \ccM,
   \]
   where $N_{Z\subset V}$ denotes the normal sheaf of $Z$ inside $V$.
   By \ref{item_Pic_of_total_space_of_bundle} again 
     we must have $\pi^*\ccM = s(B)-s'(B)$.
\end{prf}

\begin{lemma}\label{lem_contracted_section_pullback_unique}
   In Setting~\ref{set_bundle_with_contraction} 
   \begin{enumerate}
    \item \label{item_averaging_for_contraction}  
      there exists a section $s'\colon B \to U$ 
      disjoint from $s(B)$, and
    \item \label{item_pullbacks_of_line_bundles}
    the image of $\xi^*\colon \Pic Y \to \Pic U$
      is contained in $\ZZ \cdot s'(B)$.
   \end{enumerate}
   In particular, for a line bundle $\ccL$ on $Y$ 
     the linear equivalence of $\xi^*\ccL$ 
     is uniquely determined by its intersection number with a general fibre of $\pi$.
\end{lemma}

\begin{prf}
   Suppose $Z\subset Y$ is a very ample Cartier divisor in $Y$ which does not contain $y = \xi(s(B))$.
  Then $\xi^* Z$ is an effective 
    non-zero Cartier divisor on $U$ which avoids $s(B)$.
  By \cite[Lem.~A8.5]{wisniewski_length_of_extremal_rays_and_adjunction}
  there exists a section $s'\colon B \to U$ disjoint from $s(B)$, as claimed in \ref{item_averaging_for_contraction}.
  
  Let $\ccL$ be a line bundle on $Y$.
  In \ref{item_pullbacks_of_line_bundles} we claim 
    that $\xi^* \ccL = k \cdot s'(B)$ 
    for some integer $k$.
  By Lemma~\ref{lem_P1_bundle_with_two_sections}\ref{item_Pic_of_total_space_of_bundle}
  this claim is equivalent to $s^*\xi^* \ccL \simeq \ccO_B$.
  The latter is true for very ample line bundles $\ccL$
   by the argument above applied to $Z$ in the linear system of $\ccL$.
  Any ample line bundle $\ccL$ 
     is a difference $(m+1)\ccL  - m \ccL$
     with both $m\ccL$ and $(m+1) \ccL$ very ample for any sufficiently large $m$ \cite[Thm~1.2.6(i),(iv)]{lazarsfeld}.
  Therefore the claim also holds for ample line bundles.
  Finally, since $Y$ is projective, 
    any line bundle can be written as a difference of two ample line bundles, which concludes the proof.
\end{prf}

\begin{lemma}\label{lem_preimage_of_y_is_sB}
  In Setting~\ref{set_bundle_with_contraction} 
    suppose in addition that for every fibre $\PP^1 \subset U$ the restriction $\xi|_{\PP^1}$ is immersive and injective.
  Then the scheme-theoretic preimage $\xi^{-1}(y)$ is $s(B)$.
  In particular, the preimage is reduced and a Cartier divisor.
\end{lemma}
\begin{prf}
   The inclusion $s(B) \subset \xi^{-1}(y)$ is clear, and so is the reverse inclusion for the reduced structures $\reduced{(\xi^{-1}(y))} \subset s(B)$.
   By Lemma~\ref{lem_scheme_structures_are_equal}  
      it is enough to prove that the intersection 
      $\xi^{-1}(y)\cap \pi^{-1}(b)$ is reduced for any $b\in B$.
   We have:
   \[
      \xi^{-1}(y)\cap \pi^{-1}(b) = \left(\xi|_{\pi^{-1}(b)}\right)^{-1}(y),
   \]
   where the restriction $\xi|_{\pi^{-1}(b)}$ is immersive and injective.
   Thus the intersection $\xi^{-1}(y)\cap \pi^{-1}(b)$ is indeed reduced since $\set{y}$ is reduced.
\end{prf}

From now on we will assume in addition to
Setting~\ref{set_bundle_with_contraction} 
that $y$ is a smooth point of $Y$ 
  and the map $\xi$ lifts to a closed embedding 
  $\tilde{\xi} \colon U \to \tilde{Y}$, 
  where $\tilde{Y} \to Y$ be the blow-up of $Y$ at $y$.
We also need the case of reducible base for the applications 
  to the proof of Theorem~\ref{thm-properties-of-Cx}. 

\begin{setting}\label{set_reducible_base_immersive_into_blup}
   Suppose $B=B_1 \sqcup \dotsb \sqcup B_k$ 
    is a projective normal reduced scheme with irreducible components 
    $B_i$, $Y$ is a normal projective variety, 
    $\pi \colon U \to B$ is an analytically locally trivial $\PP^1$-bundle with a section $s\colon B \to U$ and a map $\xi\colon U \to Y$ such that $s(B)$ is contracted to a smooth point $y\in Y$.
   In particular, restricting to each component $B_i$ we obtain the situation of Setting~\ref{set_bundle_with_contraction}.
   Let $\tilde{Y} \to Y$ be the blow-up of $Y$ at $y$.
   We assume in addition that:
   \begin{itemize}
    \item $\xi$ lifts to a regular immersive and injective map
    $\tilde{\xi}\colon U \to \tilde{Y}$, and 
    \item each fibre $\PP^1 \subset U$ 
          is immersively mapped to $Y$ via  $\xi$.
    \end{itemize}
\end{setting}

\begin{lemma}\label{lem_tangent_cone_is_B}
  In Setting~\ref{set_reducible_base_immersive_into_blup}
      the projectivised tangent cone 
      $\PP\tau_y (\xi(U))$ is isomorphic to $B$ via the map 
      $\tilde{\xi}\circ s$.
      In particular, the projectivised tangent cone is reduced.
\end{lemma}

\begin{prf}
   The projectivised tangent cone $\PP\tau_y (\xi(U))$
   is equal (as a scheme) to the intersection of $\tilde{\xi}(U)$ and the exceptional divisor of the blow-up \cite[p.~254]{harris}.
   Since $\tilde{\xi}\colon U \to \tilde{Y}$ is a closed embedding,
   this intersection is isomorphic via $\tilde{\xi}$ to the preimage under 
   $\tilde{\xi}$ of the exceptional divisor. The latter is in turn equal to the (scheme-theoretic) preimage $\xi^{-1}(y)$. 
   Furthermore,  $\xi^{-1}(y)=s(B)$ 
      by Lemma~\ref{lem_preimage_of_y_is_sB}
      and clearly $B$ and $s(B)$ are isomorphic via $s$, 
   which concludes the proof of the lemma.
\end{prf}

For subvarieties $Z_1, Z_2 \subset \PP^N$ recall that their \emph{join} $Z_1*Z_2$
is the closure of the locus of lines between points $z_1\in Z_1 $ and $z_2 \in Z_2$.
A special case of join that is relevant to this subsection is the projective cone:
Suppose $Z \subset \PP^{N-1}$ is a projective variety, 
then by  
  $\operatorname{cone}(Z) \subset \PP^N$ we denote 
  the projective cone $Z*\set{v}$ over $Z$,
  where $v\in \PP^N \setminus \PP^{N-1}$.

\begin{prop}\label{prop_normalisations_of_cones_are_equal}
  In Setting~\ref{set_reducible_base_immersive_into_blup}
      the normalisations of 
      $ \operatorname{cone}\left(\PP\tau_y (\xi(U))\right)$
      and of $\xi(U)$ are isomorphic
        in such a way that for each $b\in B$ 
        the underlying birational map
        $ \operatorname{cone}\left(\PP\tau_y (\xi(U))\right)\dashrightarrow \xi(U)$
      takes the line in 
      $ \operatorname{cone}\left(\PP\tau_y (\xi(U))\right)$
      joining the vertex $v$ and $\tilde{\xi}\circ s(b)\in \PP\tau_y (\xi(U)$ 
      onto $\xi(\pi^{-1}(b))$.
\end{prop}

\begin{prf}
   Since $B$ is normal, it is a disjoint union of its irreducible components.
   Thus it is enough to argue for each component separately, 
      and assume that $B$ is irreducible.
   We use Lemma~\ref{lem_tangent_cone_is_B}
   to identify $\PP\tau_y (\xi(U))$ and $B$.
   The embedding $\PP\tau_y (\xi(U)) \subset \PP (T_y Y)$
     determines  a natural very ample line bundle $\ccO_{\PP\tau_y (\xi(U))}(1)$.
   We begin the proof by recognising this line bundle 
     as a line bundle on $B$. 
   Denote by $E \subset \tilde{Y}$ the exceptional divisor (preimage of $y$ under the blow up map).
   Note that
   \[
        \ccO_{\PP\tau_y (\xi(U))}(-1) =
        E|_{\PP\tau_y (\xi(U))}   
        \text{ and }
        s^* \circ \tilde{\xi}^* E = 
        s^*(s(B)) =  \ccM
   \]
   in the notation of Lemma~\ref{lem_P1_bundle_with_two_sections}.
   Thus $\ccO_{\PP\tau_y (\xi(U))}(1) = \ccM^*$ 
     via the isomorphism $\tilde{\xi}\circ s$
     and in particular, $\ccM^*$ is very ample.
     
   Therefore the normalisation of 
   $ \operatorname{cone}\left(\PP\tau_y (\xi(U))\right)$  
   is equal to $\mathfrak{C}(B, \ccM^*)$,
   the normal cone on $(B, \ccM^*)$ in the language of 
   \cite[\S1.1.8]{beltrametti_sommese_adjunction_theory_of_projective_varieties}.
   By construction,
   \begin{align*}
      \mathfrak{C}(B, \ccM^*) &
      = \Proj\left(\bigoplus_{i=0}^{\infty} H^0\left(\operatorname{Sym}^i\left(\ccO_B \oplus \ccM^*\right)\right)\right)\\
      &= \Proj\left(\bigoplus_{i=0}^{\infty} H^0\left(\operatorname{Sym}^{ik}\left(\ccO_B \oplus \ccM^*\right)\right)\right)
   \end{align*}
   (the latter equality holds for any integer $k>0$ by taking the Veronese subalgebra).
   
   Let $\ccL$ be a very ample line bundle on $Y$ 
      and suppose that 
      for a general fibre $\PP^1 \subset U$ 
      the image $\xi(\PP^1)$ has degree $k>0$ 
      with respect to $\ccL$.
   Recall from Lemma~\ref{lem_contracted_section_pullback_unique}
   that $\xi^*(i\cdot \ccL) = ik \cdot s'(B)$ and 
   from Lemma~\ref{lem_projectivisation_of_vb} 
   that $\pi_*(ik \cdot s'(B)) = \operatorname{Sym}^{ik}(\ccO_B\otimes \ccM^*)$.
   Therefore we obtain a map of algebras:
   \[
     \bigoplus_{i=0}^{\infty} H^0(Y,i \cdot \ccL) 
     \stackrel{\xi^*}{\longrightarrow}
     \bigoplus_{i=0}^{\infty} H^0(U, ik \cdot s'(B)) 
   \  \stackrel{\pi_*}{=} \
     \bigoplus_{i=0}^{\infty} H^0(B, \operatorname{Sym}^{ik}(\ccO_B\otimes \ccM^*)).
   \]
   The map induces a morphism of projective spectra:
   \[
      \zeta\colon \mathfrak{C}(B, \ccM^*) \to Y
   \]
   that factorises $\xi$ through 
   $U \to \mathfrak{C}(B, \ccM^*) \stackrel{\zeta}{\longrightarrow} Y$.
   This is because 
   \[
    U = \SheafyProj_{B} \left(\bigoplus_{i=0}^{\infty}
   \operatorname{Sym}^{ik}(\ccO_B\otimes \ccM^*)\right),
   \]
   that is, $U$ is the relative $\SheafyProj$ of a sheaf of algebras, while $\mathfrak{C}(B, \ccM^*)$ is the non-relative $\Proj$ of the sections of the same sheaf 
   (moreover,  $\ccO_B\otimes \ccM^*$ and its symmetric powers are base point free, since $\ccM^*$ is very ample).
   Furthermore, the map $U \to \mathfrak{C}(B, \ccM^*)$ 
     contracts $s(B) = \PP(\ccO_B) \subset \PP(\ccO_B\oplus\ccM) = U$ to a single point by construction. 
   
   Therefore $\zeta$ is bijective onto $\xi(U)$,
     in particular, finite and birational,
     hence $\zeta$ must be the normalisation of $\xi(U)$.
   Thus the two varieties,
   $\xi(U)$ and $ \operatorname{cone}\left(\PP\tau_y (\xi(U))\right)$, indeed have isomorphic normalisations,
      and the desired correspondence between lines follows from the naturality of the construction.
\end{prf}

Note that we do not know if the birational map
$ \operatorname{cone}\left(\PP\tau_y (\xi(U))\right)\dashrightarrow \xi(U)$ 
or its inverse is regular. 
In \cite[\S6.1]{kebekus_lines2} Kebekus claims (without proof) that under additional assumptions the map is an isomorphism. 
Unfortunately, it is not clear how to prove the claim, 
and this gap invalidates the proof of the claim about cone in \cite[Thm~1.1]{kebekus_lines2}.

\section{Loci swept out by lines}
\label{section-loci-swept-by-lines}

A rational curve $l \subset X$ is a \emph{contact line} (or simply a \emph{line}) 
if \mbox{$\deg L|_{l} =1$}.

Let $RatCurves^n(X)$ be the normalised scheme parametrising rational curves on $X$,
as in \cite[II.2.11]{kollar_book_rational_curves}.
Let $Lines(X) \subset RatCurves^n(X)$ be the subscheme parametrising lines.
Then every component of $Lines(X)$ is a minimal component of $X$ in the sense 
of \cite{hwang_mok_birationality}.
We fix $\ccH \ne \emptyset$ a union of some irreducible components of $Lines(X)$.

By a slight abuse of notation, from now on we say $l$ is a \emph{(contact) line}
if and only if $l \in \ccH$.
For simplicity, the reader may choose to restrict his attention to one of the extreme cases:
either to the case $\ccH = Lines(X)$ 
(and thus be consistent with \cite{wisniewski} and the first sentence of this section) 
or to the case where $\ccH$ is one of the irreducible components of $Lines(X)$ 
(and thus be consistent with \cite{kebekus_lines1, kebekus_lines2}).
In general it is expected that $Lines(X)$ 
(with $X$ as in Notation~\ref{notation-with-assumptions-about-contact-mflds})
is irreducible and all the cases are the same.

\subsection{Legendrian varieties swept by lines}\label{section-review-lines}
 
We denote by $C_x \subset X$ the locus of contact lines through $x \in X$.
Let $\legPF{}_x:= \PP\tau_x C_x \subset \PP(T X)$.
Note that with our assumptions both $C_x$ and $\legPF{}_x$ are closed subsets of
 $X$ or $\PP(T_x X)$ respectively.

The following theorem briefly summarises and corrects results of \cite{kebekus_lines2} and earlier:

\begin{thm} \label{thm-properties-of-Cx}
  With $X$ as in Notation~\ref{notation-with-assumptions-about-contact-mflds} 
  let $x\in X$ be any point.
  Then:
  \begin{enumerate}
     \item\label{item-exist-lines}
        There exist lines through $x$, in particular $C_x$ and $\legPF{}_x$ are non-empty.
     \item\label{item-Cx-Yx-Legendrian}  
        $C_x$ is Legendrian in $X$ and $\legPF{}_x \subset \PP(F_x)$ and $\legPF{}_x$
        is Legendrian in $\PP(F_x)$.
     \item\label{item-Yx}
        If in addition $x$ is a general point of $X$,
        then $\legPF{}_x$ is smooth.
        Further, all lines through $x$ are smooth and two different lines intersecting at $x$ do 
        not intersect anywhere else, nor they do share a tangent direction.
        Every point in $\legPF{}_x$ is a tangent direction to a contact line through $x$.
        Moreover the normalisations 
        of $\operatorname{cone}(\legPF{}_x)$ and $C_x$ 
        are isomorphic in such a way that
        the natural birational map 
        $\operatorname{cone}(\legPF{}_x)\dashrightarrow C_x$ 
        maps lines through the vertex of the cone 
        to contact lines through $x$.
\end{enumerate}
\end{thm}

\begin{prf}
   Part~\ref{item-exist-lines} is proved in 
      \cite[\S2.3]{kebekus_lines1}.

   The proof of \ref{item-Cx-Yx-Legendrian} is essentially contained in \cite[Prop.~2.9]{4authors}.
   Explicit statements are in \cite[Prop.~4.1]{kebekus_lines1} for $C_x$
   and in \cite[Lemma~5]{wisniewski} for $\legPF{}_x$.
   Also \cite{hwang_mok_VMRTs} may claim the authorship of this observation,
   since the proof in the homogeneous case is no different than in the general case.
   
   Assume  $x\in X$ is a general point.
   The statements of \ref{item-Yx} are basically \cite[Thm~1.1]{kebekus_lines2},
   which however assumes (in the statement) that $\ccH$ is irreducible.
   This is never used in the proof, with the exception of the argument for the irreducibility of $C_x$
   --- see however 
   Remark~\ref{remark-kebekus-irreducibility}.
   
   More explicitly, 
   let $B$ be the scheme parametrising contact lines through $x$.
   Since every curve through $x$ is free \cite[Lem.~3.5]{kebekus_lines1}, 
      thus $B$ is smooth by 
      \cite[Thms~II.1.7, II.2.16]{kollar_book_rational_curves}. 
   (See also \cite[Prop.~7.3]{jabu_kapustka_kapustka_special_lines}
     for a related discussion).
   We are going to verify that we are in the situation of 
   Setting~\ref{set_reducible_base_immersive_into_blup}
     with $Y=X$, $y=x$, $B$ as above, 
     $U$ the universal space of lines through $x$,
     $\xi\colon U \to X$ is the evaluation map,
     and $C_x=\xi(U)$. 
   Indeed, $B$ is smooth, hence normal and reduced,
     the section $s\colon B \to U$ assigns to $b\in B$ the point $x$ in the corresponding $\PP^1\subset X$ (in particular, $\xi$ contracts the section to $x$).
   Every line through $x$ is smooth \cite[Prop~3.3]{kebekus_lines1},
     hence $\xi^{-1}(x)= s(B)$ by Lemma~\ref{lem_preimage_of_y_is_sB}.
   In particular, $\xi^{-1}(x)$ is a Cartier divisor, 
     and by  
     \cite[Prop.~II.7.14]{hartshorne} 
     (the universal property of blowing up)
     $\xi$ factorises through $\tilde{\xi} \colon U \to \tilde{X}$, where $\beta\colon \tilde{X}\to X$ is the blow-up of $X$ at $x$.
   The only thing left to check from
   Setting~\ref{set_reducible_base_immersive_into_blup}
   is that $\tilde{\xi}$ is injective and immersive.
   
   To check the injectivity we take two distinct points 
      $u_1, u_2 \in U$
   and suppose $\tilde{\xi}(u_1)= \tilde{\xi}(u_2)$.
   In particular, $\xi(u_1)= \xi(u_2)$.
   Since $\xi$ is injective upon restricting to any fibre, 
     we must have $\pi(u_1) \ne \pi(u_2)$, 
     and either both $u_1, u_2\notin s(B)$
     or both $u_1, u_2\in s(B)$.
   By \cite[Thm~5.1]{kebekus_lines2} the former case is impossible,
     and 
     by \cite[Thm~4.1]{kebekus_lines2} the latter case is impossible,
     concluding the proof of injectivity of $\tilde{\xi}$.
  
  To check that $\tilde{\xi}$ is immersive, 
     we note that $\pi\colon U\to B$ is flat, and since 
     $\tilde{\xi}|_{\pi^{-1}(b)}$ is an embedding for any fibre 
     $\pi^{-1}(b)$, 
     thus $\tilde{\xi} \times \pi\colon U \to \tilde{X} \times B$
     is an embedding. 
  Therefore, by the universal property of Hilbert scheme,
     there is a natural map $B \to \Hilb(\tilde{X})$.
  This map is injective, since each $b$ represents a different curve in $\tilde{X}$.
  The tangent space to $\Hilb(\tilde{X})$ at the image of $b$
     is 
     $H^0\left(
     N (\PP^1 \subset \tilde{X})
     \right)$, where $\PP^1$ is $\tilde{\xi}(\pi^{-1}(b))$ 
     and $N$ denotes the normal bundle 
     \cite[Thm I.2.8.1]{kollar_book_rational_curves}.
  We will now calculate this tangent space.
  Let $E\subset \tilde{X}$ denote the exceptional divisor of the blow-up $\beta$. Then we have the following short exact sequence \cite[Lem.~15.4(iv)]{fulton_intersection_theory}:
  \[
     0 \to  T \tilde{X} \to \beta^*TX  \to T E (E) \to 0.
  \]
  Then we can restrict this sequence to 
  $\PP^1=\tilde{\xi}(\pi^{-1}(b))$ and obtain:
  \[
     0 \to  T \tilde{X}|_{\PP^1} \to \ccO_{\PP^1}(2,1^{n-1},0^{n+1})  \to T_{\tilde{\xi}(s(b))} E \to 0
  \]
  where the middle term is $\ccO_{\PP^1}(2,1^{n-1},0^{n+1}) = \ccO(2)\oplus \ccO(1)^{\oplus (n-1)}\oplus \ccO^{\oplus (n+1)}$ 
  by \cite[Lem.~3.5]{kebekus_lines1} 
  and  $T_{\tilde{\xi}(s(b))} E$ is the skyscraper sheaf supported at $\xi(s(b))$ (the intersection point of $\PP^1$ and $E$) 
    with the fibre $T_{\tilde{\xi}(s(b))} E$.
  Since $\PP^1$ is transversal to $E$
    and the $\ccO(2)$ corresponds to the tangent space to $\PP^1$, 
    it follows that 
    $T \tilde{X}|_{\PP^1}= \ccO_{\PP^1}(2,0^{n-1},(-1)^{n+1})$. 
  Therefore, $N (\PP^1 \subset \tilde{X}) = \ccO_{\PP^1}(0^{n-1},(-1)^{n+1})$ and  
  \[
    \dim H^0\left(
     N (\PP^1 \subset \tilde{X})
     \right) = n-1 = \dim B.
  \]
  We conclude that the connected component of 
  the Hilbert scheme containing the image of $B$ is smooth and isomorphic to $B$ and $U$ is the restriction of the universal subscheme to this component.
  By the decomposition of the normal bundle its sections do not vanish anywhere.
  Thus $\tilde{\xi}$ is immersive by 
  \cite[Prop.~II.3.4]{kollar_book_rational_curves}, 
  taking into account the correspondence between 
  the $\Hom$-scheme and the Chow variety 
  \cite[Thm II.2.16]{kollar_book_rational_curves} and the fact that the Hilbert to Chow is an isomorphism on this component 
  \cite[Thms I.3.17, I.3.21]{kollar_book_rational_curves}.
  
Once we have estalished that this is the situation of 
Setting~\ref{set_reducible_base_immersive_into_blup}, 
the claim of \ref{item-Yx} follows from 
Lemma~\ref{lem_tangent_cone_is_B}
and Proposition~\ref{prop_normalisations_of_cones_are_equal}.  
\end{prf}

\begin{rem}\label{remark-kebekus-irreducibility}
Note that (assuming $\ccH$ is irreducible)
Kebekus \cite{kebekus_lines2} also stated that  $C_x$ and $\legPF{}_x$ 
are irreducible for general $x$.
However it was observed by Kebekus himself together with the author that 
there is a gap in the proof.
This gap is on page 234 in Step~2 of proof of Proposition~3.2 
where Kebekus claims to construct ``a well defined family of cycles'' parametrised by a divisor $D^0$.
This is not necessarily a well defined family of cycles:
Condition (3.10.4) in \cite[\S{}I.3.10]{kollar_book_rational_curves}
is not necessarily satisfied if $D^0$ is not normal 
and there seem to be no reason to expect that $D^0$ is normal.
As a consequence the map $\Phi\colon D^0 \to \operatorname{Chow}(X)$
is not necessarily regular at non-normal points of $D^0$ and it might contract some curves.
Another problem in \cite{kebekus_lines2} concerns the claim that $C_x$ is a cone over $\legPF{}_x$ --- we discuss this problem in detail in \S\ref{sec_bundles_of_P1s}.
\end{rem}

Let us define:
\begin{align*}
 C^2 & \subset X \times X \\
 C^2 & := \set{ (x, y) \mid  y \in C_{x}},
\end{align*}
i.e.~this is the locus of those pairs $(x, y)$, which are both on the same contact line.
Again this locus is a closed subset of $X \times X$.

Analogously, define:
\[
   C^3 : = C^2 \times_X C^2
\]
so that:
\begin{align*}
 C^3 & \subset X \times X  \times X\\
 C^3 & := \set{ (x, y, z) \mid  y \in C_{x}, z \in C_{y}}.
\end{align*}
Finally, for $x \in X$ we also define $C^2_x$:
\begin{align*}
 C^2_x & \subset X  \times X  \simeq \set{x} \times X \times X \\
 C^2_x & := \set{ (y, z) \mid  y \in C_{x}, z \in C_{y}},
\end{align*}
with the scheme structure of the fibre of $C^3$ under the projection on the first coordinate.
Since for all $x \in X$ all irreducible components of $C_x$ are of dimension~$n$
 (see Theorem~\ref{thm-properties-of-Cx})
we conclude:

 \begin{prop}\label{prop-dimensions-of-Cs}
  All $C^2$, $C^2_x$, $C^3$ are projective subschemes, they are all of pure dimension,
  and their dimensions are:
  \begin{itemize}
   \item $\dim C^2 = 3n+1$.
   \item $\dim C^2_x = 2n$.
   \item $\dim C^3 = 4n+1$.
  \end{itemize}
\end{prop}
\noprf

\subsection{Joins and secants of Legendrian subvarieties}

Let $Y_1, Y_2 \subset \PP^N$ be two subvarieties. 
Note that the expected dimension of the join $Y_1*Y_2$ is $\dim Y_1 + \dim Y_2 +1$.
We are only concerned with two special cases: either $Y_1$ and $Y_2$ are disjoint or
$Y_1=Y_2$.

\begin{lemma}\label{lemma-join-of-two-disjoint-varieties}
  If $Y_1, Y_2 \subset \PP^N$ are two disjoint subvarieties of dimensions $k-1$ and $N-k$
  respectively, then their join $Y_1 * Y_2$ fills out the ambient space,
  i.e.~this join is of expected dimension.
\end{lemma}

\begin{prf}
  Let $p\in \PP^N$ be a general point and consider the projection 
  $\pi:\PP^N \dashrightarrow \PP^{N-1}$ away from $p$.
  Let $Z_i = \pi(Y_i)$ for $i =1,2$.
  Since $p$ is general, $\dim Z_i = \dim Y_i$ and thus $Z_1 \cap Z_2$ is non-empty.
  Let $q\in Z_1 \cap Z_2$ be any point.
  The preimage $\pi^{-1}(q)$ is a line in $\PP^N$ intersecting both $Y_1$ and $Y_2$ and passing through $p$.
\end{prf}

Recall, that the special case of join is when $Y=Y_1 =Y_2$
and $\sigma_2(Y) := Y*Y$ is the \emph{secant variety} of $Y$.

\begin{prop}\label{prop-secant-of-Legendrian}
     \begin{itemize}
      \item Let $Y \subset \PP^{2n-1}$ be an irreducible linearly non-degenerate Legendrian variety.
            Then $\sigma_2(Y) = \PP^{2n-1}$.
      \item Let $Y_1,Y_2 \subset \PP^{2n-1}$ be two disjoint Legendrian subvarieties.
            Then $Y_1*Y_2 = \PP^{2n-1}$.
     \end{itemize}
\end{prop}

\begin{prf}
     If $Y$ is irreducible, then this is proved in the course of proof of Prop.~17(2) 
     in \cite{landsbergmanivel04}.
     
     If $Y_1$ and $Y_2$ are disjoint,
     then the result follows from Lemma~\ref{lemma-join-of-two-disjoint-varieties}.
\end{prf}

\subsection{Divisors swept by broken lines}\label{section-define-Dx}

Following the idea of \Wisniewski{} \cite{wisniewski} 
we introduce the locus of broken lines 
(or reducible conics, or chains of $2$ lines) through $x$:
\[
 {D}_x := \bigcup_{y \in C_x} C_y.
\]
Note that ${D}_x$ is a closed subset of $X$ 
as it can be interpreted as the image of projective variety $C^2_x \subset X \times X$ under a proper map, 
which is the projection onto the last coordinate. 
By analogy to the case of lines consider also:
\begin{align*}
 {D}^2 & \subset  X \times X \\
 {D}^2 & := \set{ (x, z) \mid   \exists_{y \in C_{x}} \text{ s.t. } z \in C_{y}},
\end{align*}
i.e.~${D}^2$ is the projection of $C^3$ onto first and third coordinates.
Thus again ${D}^2$ is a closed subset of the product. 
Set theoretically ${D}_x$ is the fibre over $x$ of
(either of) the projection ${D}^2 \to X$ and if we consider ${D}^2$ as a reduced scheme,
then we can assign to ${D}_x$ the scheme structure of the fibre.

It follows immediately from the above discussion and Proposition \ref{prop-dimensions-of-Cs},
that every component of ${D}_x$ has dimension at most $2n$ and 
every component of ${D}^2$ has dimension at most $4n+1$.
In fact the equality holds.

\begin{thm}
  \label{thm-Dx-is-a-divisor}
  Let $x \in X$ be any point. 
  Then the locus ${D}_x$ is of pure codimension~$1$.
\end{thm}

\begin{prf}
   Assume first that $x\in X$ is a general point.
   Recall, that $C^2_x \subset X \times X$ has two projections:
   \[
     \xymatrix{ C^2_x  \ar@{->>}[r]^{\pi_2} \ar@{->>}[d]^{\pi_1} & D_x \\
                C_x}
   \]
   Fix $\component{D_x}$ to be an irreducible component of $D_x$.
   Then $\component{D_x}$ is dominated 
   by some component $\component{C^2_x}$ of $C^2_x$.
   Dimension of $\component{C^2_x}$ is equal to $2n$ by Proposition~\ref{prop-dimensions-of-Cs}.

   For $y \in C_x$ the fiber ${\pi_1}^{-1}(y) \subset C^2_x$
   is equal to $\set{y} \times C_y$.
   In particular, 
   by Theorem~\ref{thm-properties-of-Cx}\ref{item-Cx-Yx-Legendrian}
   the fibers of $\pi_1$ have constant dimension $n$. 
   Thus $\component{C^2_x}$ is mapped onto an irreducible component $\component{C_x}$ of $C_x$.
   Finally, let $C'$ be an irreducible component of the preimage ${\pi_1}^{-1}(x)$ 
   which is contained in $\component{C^2_x}$.
   Note that $C'$ can be identified with an irreducible component of $C_x$,
   because ${\pi_1}^{-1}(x) = \set{x} \times C_x$.
  
   We claim that the projectivised tangent cone $\PP \tau_x \component{D_x}$ contains 
   the join of two tangent cones 
   \[
    (\PP\tau_x C') * (\PP\tau_x \component{C_x}) \subset \PP F_x \subset \PP T_x X.
   \]
   The proof of the claim is a baby version of \cite[Thm~3.11]{hwang_kebekus_chains}.
   There however Hwang and Kebekus assume $C_x$ is irreducible and thus their 
   results do not neccessarily apply directly here.
   Let $l_0$ be a general line through $x$ contained in $C'$ 
   and let $l$ be a general line through $x$ contained in $\component{C_x}$.
   To prove the claim it is enough to show there exists 
   a surface $S\subset D_x$ containing $l_0$ and $l$ which is smooth at $x$, 
   since in such a case $T_x S \subset \tau_x D_x$ and  
   $\PP T_x S$ is the line between $\PP T_x l$ and $\PP T_x l_0$.
   
   We obtain $S$ by varying $l_0$.
   Consider $\ccH_l \subset \ccH$ the parameter space for lines on $X$, 
   which intersect $l$. 
   By Theorem~\ref{thm-properties-of-Cx}\ref{item-Yx} the space $\ccH_l$ comes with a projection 
   $\xi\colon \ccH_l \dashrightarrow l$, 
   which maps $l'\in \ccH_l$ to the intersection point of $l$ and $l'$,
   and which is well defined on an open subset containg all lines through $x$.

   By generality of our choices, $l_0$ is a smooth point of $\ccH_l$ and $\xi$ is submersive at $l_0$.
   In the neighbourhood of $l_0$ choose a curve $A \subset \ccH_l$ smooth at $l_0$ 
   for which $\xi|_A$ is submersive at $l_0$.
   Then the locus in $X$ of lines which are in $A$ sweeps a surface $S \subset X$,
   which is smooth at $x$, contains $l_0$, and contains an open subset of $l$ around $x$.
   Thus the claim is proved and:
   \begin{equation}\label{equ-join-is-contained-in-tangent-to-Dx}
      (\PP\tau_x C') * (\PP\tau_x \component{C_x}) \subset \PP \tau_x \component{D_x}
   \end{equation}

   Now we claim that $F_x \subset \tau_x D_x$. 
   For this purpose we analyse three cases. 

   In the first case $C' = \component{C_x}$ and 
   the corresponding component of $\legPF{}_x$ is nondegenerate in $\PP F_x$.
   Then $\PP\tau_x C'$ is non-degenerate in $\PP F_x$ by 
    Theorem~\ref{thm-properties-of-Cx} 
   and thus
   \[
     (\PP\tau_x C') * (\PP\tau_x \component{C_x}) = \sigma_2(\PP\tau_x C') = \PP(F_x)
   \]
   by 
    Proposition~\ref{prop-secant-of-Legendrian}. 
   Combining with 
    \eqref{equ-join-is-contained-in-tangent-to-Dx}
   we obtain the claim.
   
   In the second case $C'$ and $\component{C_x}$ are different components of $C_x$.
   Then by generality of $x$ and by 
    Theorem~\ref{thm-properties-of-Cx},
   the two tangent cones $(\PP\tau_x C')$ and $(\PP\tau_x \component{C_x})$ are disjoint.
   Thus again 
   \[
     (\PP\tau_x C') * (\PP\tau_x \component{C_x}) = \PP(F_x)
   \]
   by 
    Proposition~\ref{prop-secant-of-Legendrian}. 
   Combining with 
    \eqref{equ-join-is-contained-in-tangent-to-Dx}
   we obtain the claim.

Finally, we exclude the possibility that 
   $C' = \component{C_x}$ and 
   the corresponding component of $\legPF{}_x$
   is contained in a hyperplane $H\subset \PP F_x$.
Then, since the component of $\legPF{}_x$ 
   is smooth by Theorem~\ref{thm-properties-of-Cx}\ref{item-Yx}, 
   it must be a linear space by 
   \cite[Thm~3.4]{jabu06} or \cite[Prop.~17, item~1]{landsbergmanivel04}.
By Theorem~\ref{thm-properties-of-Cx}\ref{item-Yx} the normalisation of $ \component{C_x}$ is $\iota \colon \PP^n \to \component{C_x}$ and since contact lines through $x$ in $ \component{C_x}$ are pulled back to ordinary lines in $\PP^n$, 
  we must have $\iota^*L|_{\component{C_x}} = \ccO_{\PP^n}(1)$.
Therefore images of all lines in $\PP^n$ are also contact lines in $X$ and since $C' = \component{C_x}$,
  we would have $\component{D_x} = \component{C_x}$.
However, there are other lines through $x$ by \cite[Thm~4.4]{kebekus_lines1} and $\component{C_x}$ is strictly contained in another component of $D_x$, which is a contradiction, as it shows that $\component{D_x}$ is not an irreducible component of $D_x$.

   Thus in any case for a general $x \in X$, every component of $D_x$ has dimension at least $2n$.
   The dimension can only jump up
   at special points when one has a fibration,
   thus also at special points every component of $D_x$ has dimension at least $2n$.
   Earlier we observed that $\dim D_x \le 2n$, thus the theorem is proved.
\end{prf}

\begin{prop}\label{prop-what-is-Dx-for-X-homogeneous}
  If $X$ is the adjoint variety of $G$, and $x\in X$, then $D_x$ 
  is the hyperplane section of $X \subset \PP(\gotg)$
  perpendicular to $x$ via the Killing form.
\end{prop}

\begin{prf}
  Let $X = G/P$, where $P$ is the parabolic subgroup preserving $x$.
  Notice, that $D_x$ must be reduced (because $D$ is reduced and $D_x$ is a general fibre of $D$).
  Also $D_x$ is $P$-invariant, because the set of lines is $G$ invariant 
  and $D_x$ is determined by $x$ and the geometry of lines on $X$.
  We claim, there is a unique $P$-invariant reduced divisor on $X$,
  and thus it must be the hyperplane section as in the statment of proposition.

  So let $\Delta$ be a $P$-invariant divisor linearly equivalent to $L^k$ for some $k \ge 0$.
  Also let $\rho_{\Delta}$ be a section of $L^k$ which determines $\Delta$. 
  The module of sections $H^0(L^k)$ is an irreducible $G$-module by Borel-Weil theorem 
  (see \cite{serre_representations}), with some highest weight  $\omega$.
  Since the Lie algebra $\gotp$ of $P$ contains all positive root spaces, 
  by \cite[Prop.~14.13]{fultonharris} there is a unique $1$-dimensional $\gotp$-invariant
  submodule of $H^0(L^k)$, it is the highest weight space $H^0(L^k)_{\omega}$.
  So $\rho_{\Delta} \in H^0(L^k)_{\omega}$ and $\Delta$ is unique.
  
  The hyperplane section of $X \subset \PP(\gotg)$
  perpendicular to $x$ via the Killing form is a divisor in $\lvert L \rvert$,  and it is $P$-invariant,
  and so are its multiples in $\lvert L^k \rvert$.
  So by the uniqueness $\Delta$ must be equal to $k$ times this hyperplane section.
  Thus $\Delta$ is reduced if and only $k=1$ and so $D_x$ is the hyperplane section.
\end{prf}

\subsection{Tangent bundles restricted to lines}

   Let $l$ be a  line through a general point $y\in X$.
   Recall from  \cite[Fact~2.3]{kebekus_lines2} that: 
   \begin{align}
      T X|_{l} & \simeq {\ccO_l}(2) \oplus {\ccO_l}(1)^{n-1} \oplus {\ccO_l}^{n-1} \oplus {\ccO_l}^{2} 
        \nonumber
        \label{equ-splitting-type-of-TX-on-l}\\ 
        \nonumber
        F|_{l} & \simeq {\ccO_l}(2) \oplus {\ccO_l}(1)^{n-1} \oplus {\ccO_l}^{n-1} \oplus {\ccO_l}(-1) \\
        \nonumber
      T l      & \simeq {\ccO_l}(2)\\
   \intertext{and for general $z \in l$:}
        \nonumber
      T C_z|_{l \setminus \set{z}} & \simeq {\ccO_l}(2) \oplus {\ccO_l}(1)^{n-1}.
   \end{align}
   
   If $x \in X$ is a general point and $y \in C_x$ 
    is a general point of any of the irreducible components of $C_x$
   and $l$ is a line through $y$,
   then we want to express $T D_x |_{l}$ in terms of those splittings.
   In a neighbourhood of $l$ the divisor $D_x$ is swept by deformations $l_t$ of $l =l_0$ 
   such that $l_t$ intersects $C_x$.   
   By the standard deformation theory argument taking derivative of $l_t$ by $t$ at a point $z \in l$,
   we obtain that:
   \begin{equation}\label{equ-for-tangent-space-of-Dx}
     T_{z} D_x
      \supset \set{s(z) \in T_z X  \mid  \exists s \in H^0(T X|_{l}) \  \text{ s.t. } s(y) \in T_y C_x}  
   \end{equation}
   Moreover, at a general point $z$ we have equality in \eqref{equ-for-tangent-space-of-Dx}.
   If we mod out $T X|_l$ by the rank $n$ positive bundle $(T X|_l)^{>0}:= {\ccO_l}(2) \oplus {\ccO_l}(1)^{n-1}$, 
   then we are left with a trivial bundle ${\ccO_l}^{n+1}$. 
   Thus, since by Theorem~\ref{thm-Dx-is-a-divisor} the dimension of $T_z D_x = 2n$ for general $z \in l$,
   the vector space $T_y C_x$ must be transversal to  $(TX|_l)^{>0}$ at $y$.
   In particular, if $z\ne y$,
   then dimension of the right hand side in \eqref{equ-for-tangent-space-of-Dx} is $2n$
   and thus \eqref{equ-for-tangent-space-of-Dx} is an equality for each point $z\in l$,
   such that $z$ is a smooth point of $D_x$.
   
   We conclude:
   \begin{prop}\label{prop-tangent-to-Dx}
   Let $x \in X$ be a general point and $y \in C_x$ 
   be a general point of any of the irreducible components of $C_x$
   and $l$ be any line through $y$.
   Then there exists a subbundle $\varGamma \subset TX|_l$
   such that:
   \begin{align*}
     \varGamma & = {\ccO_l}(2) \oplus {\ccO_l}(1)^{n-1} \oplus {\ccO_l}^{n},\\
     \varGamma \cap F|_{l}
                 & = {\ccO_l}(2) \oplus {\ccO_l}(1)^{n-1} \oplus {\ccO_l}^{n-1} = (F|_l)^{\ge 0}
   \end{align*}
    and if $z \in l$ is a smooth point of $D_x$, then $T_z D_x = {\varGamma}_z$.
   \end{prop}
   \noprf

\section{Duality}
\label{section-duality}

An effective divisor $\Delta$ on $X$ is an element of divisor group 
(and thus a positive integral combination of codimension $1$ subvarieties of $X$)
and also a point in the projective space $\PP(H^0 \ccO_X(\Delta))$ or a hyperplane 
in $\PP(H^0 \ccO_X(\Delta)^*)$. 
In this section we will constantly interchange these three interpretations of $\Delta$.
In order to avoid confusion we will write:
\begin{itemize}
  \item $\divisor{\Delta}$ to mean the divisor on $X$;
  \item $\point{\Delta}$ to mean the point in $\PP(H^0 \ccO_X(\Delta))$ 
                                                     or in a fixed linear subsystem.
  \item $\hyperplane{\Delta}$ to mean the hyperplane in $\PP(H^0 \ccO_X(\Delta)^*)$ 
                                                     or in dual of the fixed subsystem.
\end{itemize}

In \S\ref{section-define-Dx} we have defined $D \subset X \times X$,
which we now view as a family of divisors on $X$ parametrised by $X$.
Since the Picard group of $X$ is discrete and $X$ is smooth and connected,
it follows that all the divisors $D_x$ are linearly equivalent. 
Thus let $E \simeq L^{\otimes k}$ be the line bundle $\ccO_X(D_x)$.
Consider the following vector space $\LinSyst \subset H^0(E)$:
\[
  \LinSyst := \sspan \set{s_x : x \in X} \text { where $s_x$ is a section of $E$ vanishing on $D_x$.}  
\]
Hence $\PP \LinSyst$ is the linear system spanned by all the $D_x$.

Further, consider the map
\[
  \phi: X \to \PP \LinSystDual
\]
determined by the linear system $\LinSyst$,
i.e.~mapping point $x \in X$ to the hyperplane in $\PP \LinSyst$
consisting of all divisors containing $x$.
\begin{rem}\label{rem-phi-is-regular}
Note that $\phi$ is regular, since for every $x \in X$ there exists $w \in X$, 
such that $x \notin D_w$
(or equivalently, $w \notin D_x$).
\end{rem}
Since $E$ is ample, it must intersect every curve in $X$ and hence
$\phi$ does not contract any curve. Therefore $\phi$ is finite to one. 

\begin{prop} \label{prop-homogeneous-iff-k-eq-1}
  If $X$ is an adjoint variety, then $k=1$, i.e.~$E \simeq L$.
  If $k=1$ and the automorphism group of $X$ is reductive,
  then $X$ is isomorphic to an adjoint variety.
\end{prop}
\begin{prf}
  If $X$ is the adjoint variety of $G$, and $x\in X$, then $D_x$ 
  is the hyperplane section of $X \subset \PP(\gotg)$
  by Proposition~\ref{prop-what-is-Dx-for-X-homogeneous}.

  If $k=1$ and the automorphism group of $X$ is reductive,
  since $\phi$ is finite to one, we can apply 
  Beauville Theorem \cite{beauvillefano}.
  Thus $X$ is isomorphic to an adjoint variety.
\end{prf}

\subsection{Dual map}\label{section-dual-map}

In algebraic geometry it is standard to consider maps 
determined by linear systems (such as $\phi$ defined above).
However in our situation, we also have another map determined by
the family of divisors $D$.
Namely:
\begin{align*}
  \psi\colon X & \to \PP \LinSyst \\
             x & \mapsto \point{D_x}.
\end{align*}
So let $\ccS \subset \ccO_X \otimes \LinSystDual \simeq  X  \times \LinSystDual$
be the pullback under $\phi$ of the universal hyperplane bundle, i.e. 
the corank 1 subbundle
such that the fibre of $\ccS$ over $x$ 
is $\hyperplane{D_x} \subset \LinSystDual$.
We note that $\PP(\ccS)$  is both a projective space bundle on $X$
 and also it is a divisor on $X\times \PP \LinSystDual$.
Also $D= (\id_X \times \phi)^* \PP(\ccS)$ as divisors.

We can also consider the line bundle dual to the cokernel of $\ccS \to \ccO_X \otimes \LinSystDual$,
i.e.~the  subbundle ${\ccS}^{\perp} \subset  \ccO_X \otimes \LinSyst$. 
This line subbundle determines section $X \to X\times \PP\LinSyst$,
where $x \mapsto (x , \point{D_x})$. 
So $\psi$ is the composition of the section and the projection:
\[
  X \to X\times \PP\LinSyst \to \PP\LinSyst.
\]

Every map to a projective space is determined by some linear system. 
We claim the $\psi$ is determined by $\LinSyst$,
precisely the system that defines $\phi$
and thus that there is a natural linear isomorphism between 
$\PP \LinSyst$ and $\PP \LinSystDual$.
 
\begin{prop}
We have $\psi^* \ccO_{\PP \LinSyst}(1) \simeq E$ and 
the linear system cut out by hyperplanes 
\[
 \psi^* H^0 \left( \ccO_{\PP \LinSyst}(1) \right) := \set{\psi^* s : s \in \LinSystDual} \subset H^0(E) 
\]
 is equal to $\LinSyst$.
\end{prop}

\begin{prf}
  For fixed $x\in X$ let $\phi(x)^{\perp} \subset \PP \LinSyst$ be the hyperplane dual to 
  $\phi(x) \in \PP \LinSystDual$. 
  To prove the proposition it is enough to prove 
  \begin{equation} \label{equ-pullback-by-psi-is-Dx}
    \psi^*(\phi(x)^{\perp}) = \divisor{D_x}.
  \end{equation}
  Since we have the following symmetry property of $D$:
  \[
     x \in D_y \iff y \in D_x,
  \]
  the set theoretic version of \eqref{equ-pullback-by-psi-is-Dx}
  follows easily: 
  \[
     y \in \psi^*(\phi(x)^{\perp}) 
      \iff \psi(y) \in \phi(x)^{\perp} 
      \iff \hyperplane{D_y} \ni \phi(x)
      \iff D_y \ni x.
  \]
  However, in order to prove the equality of divisors in \eqref{equ-pullback-by-psi-is-Dx}
  we must do a bit more of gymnastics, which translates the equivalences above into 
  local equations.
  The details are below.

  The pull back of $\phi(x)^{\perp}$ by the projection $X\times \PP\LinSyst \to \PP\LinSyst$ is just 
  $X\times \phi(x)^{\perp}$.
  Then the pull-back of the product by the section $X \to X\times \PP\LinSyst$
  associated to ${\ccS}^{\perp}$ is just the subscheme of $X$ defined by 
  $\set{y \in X \mid ({\ccS}^{\perp})_y \subset  \phi(x)^{\perp}}$
  (locally, this
  is just a single equation: the spanning section of ${\ccS}^{\perp}$
  satisfies the defining equation of $\phi(x)^{\perp}$).
  But this is clearly
  equal to the dual equation $\set{y \mid \PP(\ccS_y) \ni \phi(x) }$. 
  If we  let $\rho_x$ be the section 
  \begin{align*}
    \rho_x\colon X&\to X \times X \\
    \rho_x (y) & :=  (y,x)
  \end{align*}
  then we have:
  \[
    \psi^*(\phi(x)^{\perp}) = {\rho_x}^* \circ (\id_X \times \phi)^* (\PP(\ccS)) = {\rho_x}^* (D)
    = \divisor{D_x}
  \]
  as claimed.
\end{prf}

Thus we have a canonical linear isomorphism
$f\colon \PP \LinSystDual \to \PP \LinSyst$ giving rise to the following commutative diagram:
\begin{equation} \label{equ-duality-diagram}
              \xymatrix@C1.5pc@R0.5pc{
              & & & \PP \LinSystDual \ar[dd]^{\simeq} \\
              X\ar[urrr]^{\phi} \ar[drrr]^{\psi} && \\
              & & & \PP \LinSyst.}
\end{equation}
We will denote the underlying vector space isomorphism $\LinSystDual \to \LinSyst$ 
(which is unique up to scalar) with the same letter $f$.
The choice of $f$
 combined with the canonical
pairing $\LinSyst \times \LinSystDual \to \CC$,
determines a non-degenerate bilinear form 
$B: \LinSyst \times \LinSyst \to \CC$,
with the following property:
\begin{equation}\label{equ-defining-h}
  B(\phi(x),\phi(y)) =0 \iff (x,y) \in D \iff x \in D_y \iff y \in D_x.
\end{equation}
\begin{prop}\label{prop-Killing-form-for-homogeneous}
  If $X$ is the adjoint variety of $G$,
  then $\LinSyst = H^0(L) \simeq \gotg$ and $B$ is (up to scalar) the Killing form on $\gotg$.
\end{prop}
\begin{prf}
  Follows immediately from Proposition~\ref{prop-what-is-Dx-for-X-homogeneous} 
  and Equation \ref{equ-defining-h}.
\end{prf}

\begin{cor}
          \label{cor_when_D_x_equals_D_y}
  $\phi(x) = \phi(y)$ if and only if $D_x = D_y$.
\end{cor}

\begin{prf}
  It is immediate from the definition of $\psi$
  and from Diagram \eqref{equ-duality-diagram}.
\end{prf}

\subsection{Symmetry}

Note that $B$ has the 
property that for $x\in X$,
\[
  B(\phi(x),\phi(x))=0
\]
(because $x\in D_x$).
 
\begin{prop}
              \label{h_is_sym_or_skew}
  The bilinear form $B$ is either symmetric or skew-symmetric.
\end{prop}

\begin{prf} 
  Consider two linear maps $\LinSyst \to \LinSystDual$: 
  \[
    \alpha(v):= B(v, \cdot) \ \text{ and } \  \beta(v):= B( \cdot,v).
  \]
  If $v= \phi(x)$ for some $x \in X$, then 
  \[
    \ker\bigl(\alpha(v)\bigr) 
    = \sspan\Bigl(\ker\bigl(\alpha(v)\bigr) \cap \phi(X)\Bigr)
    = \sspan\bigl(\phi(D_x)\bigr)
  \]
  and analogously $\ker(\beta(v)) = \sspan(\phi(D_x))$. 
  So  $\ker(\alpha(v)) = \ker(\beta(v))$ 
  and hence $\alpha(v)$ and $\beta(v)$ are proportional. 
  Therefore  there exists a function $\lambda: X \rightarrow \CC$ such that:
  \[
    \lambda(x) \alpha(\phi(x)) = \beta(\phi(x)).
  \]

  So for every $x,y \in X$ we have:
  \[
     B(\phi(x),\phi(y))  =
    \lambda(x) B(\phi(y),\phi(x)) = 
    \lambda(x)\lambda(y) B (\phi(x),\phi(y))
  \]
  and hence:
  \[
    \forall {(x,y) \in X\times X \setminus D}  \qquad \lambda(x)\lambda(y) = 1.
  \]
  Taking three different points we see that $\lambda$ is constant and $\lambda \equiv \pm 1$. 
  Therefore $\pm \alpha(\phi(x)) = \beta(\phi(x))$ 
  and by linearity this extends to 
  $\pm \alpha = \beta$ so $B$ is either symmetric or skew-symmetric as stated in the proposition.
\end{prf}

\begin{example}
  If $X$ is one of the adjoint varieties, then $B$ is symmetric (because the Killing form is symmetric).
\end{example}

\begin{rem}
  Consider $\PP^{2n+1}$ with a contact structure arising from a symplectic form $\omega$ on $\CC^{2n+2}$.
  Recall, that this homogeneous contact Fano manifold does not satisfy our assumptions, 
  namely, its Picard group is not generated by the equivalence class of $L$
  --- in this case $L \simeq \ccO_{\PP^{2n+1}}(2)$.
  However, \Wisniewski{} in \cite{wisniewski} considers also this generalised situation
  and defines $D_x$ to be the divisor swept by contact conics 
  (i.e.~curves $C$ with degree of $L|_C = 2$)
  tangent to the contact distribution $F$.
  Then for the projective space $D_x$ is just the hyperplane perpendicular to $x$ with respect to $\omega$.
  And thus in this case $\LinSyst = H^0\left(\ccO_{\PP^{2n+1}}(1)\right)$ 
  and the bilinear form $B$ defined from such family of divisors would be proportional to $\omega$,
  hence skew-symmetric.
\end{rem}

\begin{prf}[ of Theorem~\ref{thm-exist-Dx-and-duality}]
  $D_x$ is a divisor by Theorem~\ref{thm-Dx-is-a-divisor}.
  $\phi$ is regular by Remark~\ref{rem-phi-is-regular}.
  $\psi$ is regular by \eqref{equ-duality-diagram}.
  The non-degenerate bilinear form $B$ is constructed in \S\ref{section-dual-map}.
  It is either symmetric or skew-symmetric by Proposition~\ref{h_is_sym_or_skew}.
  In the adjoint case $B$ is the Killing form by Proposition~\ref{prop-Killing-form-for-homogeneous}.
\end{prf}

\begin{cor}
If $B$ is symmetric, then $\psi(X) \subset \PP\LinSyst$ is contained in the quadric 
 $B(v,v)=0$.
\end{cor}

\begin{cor}
If $x\in X$, then $\psi(C_x)$ is contained in a linear subspace of dimension 
at most $\left\lfloor \frac{\dim \LinSyst}{2} \right\rfloor$.
\end{cor}

\begin{prf}
If $y,z \in C_x$, then $z\in D_y$, so $B(\psi(y), \psi(z))=0$.
Therefore $\sspan(\psi(C_x))$ is an isotropic linear subspace, which
cannot have dimension bigger than $\left\lfloor\frac{\dim \LinSyst}{2}\right\rfloor$.
\end{prf}

\section{Grading}\label{section-grading}
Suppose $X \subset \PP \gotg$ is the adjoint variety of $G$. 
Assume further that a maximal torus and an order of roots in $\gotg$ has been chosen, 
then $\gotg$ has a natural grading (see \cite[\S6.1]{landsbergmanivel02}):
\[
  \gotg = \gotg_{-2} \oplus \gotg_{-1} \oplus \gotg_{0} \oplus \gotg_{1} \oplus \gotg_{2}
\]
where:
\begin{enumerate}
  \item $\gotg_{0} \oplus \gotg_{1} \oplus \gotg_{2}$ is the parabolic subalgebra $\gotp$
        of $X$.
  \item $\gotg_{0}$ is the maximal reductive subalgebra of $\gotp$.
  \item for all $i \in  \set{-2,-1,0,1,2}$ the vector space $\gotg_{i}$ is a $\gotg_0$-module.
  \item $\gotg_{2}$ is the $1$-dimensional highest root space, 
  \item $\gotg_{-2}$ is the $1$-dimensional lowest root space.
  \item The restriction of the Killing form to each
        $\gotg_{2} \oplus \gotg_{-2}$, $\gotg_{1} \oplus \gotg_{-1}$ and $\gotg_0$ 
	is non-degenerate,  
	and the Killing form $B(\gotg_i, \gotg_j)$ is identically zero 
	for $i \ne -j$.
  \item The Lie bracket on $\gotg$ respects the grading, $[\gotg_i, \gotg_j] \subset \gotg_{i+j}$
        (where $\gotg_k = 0$ for $k \notin \set{-2,-1,0,1,2}$).
\end{enumerate}

In fact the grading is determined by $\gotg_{-2}$ and $\gotg_{2}$ together with the geometry of $X$ only.
So let $X$ be as in Notation \ref{notation-with-assumptions-about-contact-mflds}
 and let $x$ and $w$ be two general points of $X$.
Define the following subspaces of $\LinSyst$:
\begin{itemize}
  \item $\LinSyst_{2}$ to be the $1$-dimensional subspace $\psi(x)$;
  \item $\LinSyst_{-2}$ to be the $1$-dimensional subspace $\psi(w)$;
  \item $\LinSyst_{1}$ to be the linear span of affine cone of $\psi(C_x \cap D_w)$;
  \item $\LinSyst_{-1}$ to be the linear span of affine cone of  $\psi(C_w \cap D_x)$;
  \item $\LinSyst_{0}$ to be the vector subspace of $\LinSyst$,
        whose projectivisation is:
        \[
          \bigcap_{y\in C_x \cup C_w} f(\hyperplane{D_y})
        \]
\end{itemize}

In the homogeneous case this is precisely the grading of $\gotg$.

\begin{prf}[ of Theorem~\ref{thm-grading-for-homogeneous}]
  First note that the classes of the 1-dimensional linear subspaces $\gotg_{2}$ and $\gotg_{-2}$ 
  are both in $X$ (as points in $\PP \gotg$).
  Moreover, they are a pair of general points, because the action of the parabolic subgroup $P < G$
  preserves $\gotg_2$ and moves freely $\gotg_{-2}$.
  This is because $\hat{T}_{[\gotg_{-2}]} X  = [\gotg_{-2} , \gotg] = [\gotg_{-2} , \gotp]$.
  
  So fix $x = [\gotg_2]$ and $w= [\gotg_{-2}]$. 
  We claim the linear span of $C_x$ (respectively $C_w$) 
  is just $\gotg_2 \oplus \gotg_1$ (respectively $\gotg_{-2} \oplus \gotg_{-1}$).
  To see that, we observe the lines on $X$ through $x$ are contained
  in the intersection of $X$ and the projective tangent space $\PP(\hat{T}_x X) \subset \PP(\gotg)$.
  In fact this intersection is equal to $C_x$:
  if $y \ne x$ is a point of the intersection,
  then the line in $\PP \gotg$ through $x$ and $y$ intersects $X$ with multiplicity at least $3$,
  but $X$ is cut out by quadrics (see for instance \cite[\S10.6.6]{procesi_book}),
  so this line must be contained in $X$. 
  Also $C_x$ is non-degenerate in $\PP(\hat{F}_x) \subset \PP(\hat{T}_x X)$.
  However $\hat{F}_x$ is a $\gotp$-invariant hyperplane in $\PP(\hat{T}_x X)$ 
  and the unique $\gotp$-invariant hyperplane in 
  \[
    \hat{T}_x X = [ \gotg, \gotg_2] = [ \gotg_{-2}, \gotg_2] \oplus \gotg_{1} \oplus \gotg_2
  \]
  is
  \[
    \hat{F}_x = [\gotg_{-1} \oplus \gotg_0 \oplus \gotg_1 \oplus \gotg_2, \gotg_2 ] = \gotg_{1} \oplus \gotg_2.
  \]
  
  Further we have seen in Proposition~\ref{prop-what-is-Dx-for-X-homogeneous}
  that $D_x$ (respectively $D_w$)
  is the intersection of  
   $\PP({\gotg_2}^{\perp_{B}}) = \PP(\gotg_2 \oplus \gotg_1 \oplus \gotg_0 \oplus \gotg_{-1})$ and $X$
  (respectively $\PP(\gotg_{-2} \oplus \gotg_{-1} \oplus \gotg_0 \oplus \gotg_{1})$ and $X$).
  Equivalently, $f(\hyperplane{D_x}) = \PP(\gotg_2 \oplus \gotg_1 \oplus \gotg_0 \oplus \gotg_{-1})$.
  Thus:
  \[
      C_x \cap D_w = C_x \cap f(\hyperplane{D_w}) 
    = C_x \cap \PP(\gotg_{-2} \oplus \gotg_{-1} \oplus \gotg_0 \oplus \gotg_{1}) 
    = C_x \cap \PP(\gotg_{1}).
  \]
  $C_x \cap \PP(\gotg_{1})$ is non-degenerate in $\PP(\gotg_{1})$,
  thus $\LinSyst_{1} = \gotg_1$ and analogously $\LinSyst_{-1} = \gotg_{-1}$.
  
  It remains to prove $\LinSyst_{0} = \gotg_0$.
  \begin{align*}
    \PP\LinSyst_{0} & = \bigcap_{y\in C_x \cup C_w} f(\hyperplane{D_z})\\
                    & = (C_x \cup C_w)^{\perp_{B}} \\
                    & = \PP(\gotg_2 \oplus \gotg_1 \oplus \gotg_{-1} \oplus \gotg_{-2})^{\perp_{B}} \\
                    & = \PP(\gotg_0).
  \end{align*}
\end{prf}

We also note the following lemma in the homogeneous case:
\begin{lemma}\label{lemma-homogeneous-then-X-cap-g1-subset-Cx}
  If $X$ is the adjoint variety of $G$,
  then
  \[
    X \cap \PP(\gotg_1) \subset C_x
  \]
  where $x$ is the point of projective space corresponding to $\gotg_2$.
\end{lemma}

\begin{prf}
  Suppose $y \in X \cap \PP \gotg_1$ 
  and let $l \subset \PP \gotg$ be the line through $x$ and $y$.
  Note that $l \subset \PP(\gotg_1 \oplus \gotg_2)$.
  Since $\gotg_1 \oplus \gotg_2 \subset [\gotg,\gotg_2] = \hat{T}_x X$, 
  hence $l \cap X$ has multiplicity at least $2$ at $x$.
  Thus $l\cap X$ has degree at least $3$ and since $X$ is cut out by quadrics,
  $l$ is contained in $X$. 
\end{prf}

\section{Cointegrable subvarieties}
\label{section-cointegrable}

\begin{defin}
  A subvariety $\Delta \subset X$ is \emph{$F$-cointegrable}
  if  $T_x \Delta \cap F_x \subset F_x$ is a coisotropic subspace for general point $x$ of each 
  irreducible component of $\Delta$.
\end{defin}
Note that this is equivalent to the definition given in \cite[\S{}E.4]{jabu_dr} 
--- this follows from the local description of the symplectic form 
    on the symplectisation of the contact manifold (see \cite[(C.15)]{jabu_dr}).

Clearly, every codimension $1$ subvariety of $X$ is $F$-cointegrable.

Assume $\Delta \subset X$ is a subvariety of pure dimension, which is $F$-cointegrable
 and let $\Delta_0$ 
be the locus where $T_x \Delta \cap F_x \subset F_x$ is a coisotropic subspace
of dimension $\dim \Delta-1$.
We define the \emph{$\Delta$-integrable distribution $\Delta^{\perp}$}
 to be the distribution defined over $\Delta_0$ by:
\[
  {\Delta^{\perp}}_x:=\left(T_x \Delta \cap F_x \right)^{\perp_{\ud \theta}} \subset F_x
\]
We say an irreducible subvariety $A \subset X$ is \emph{$\Delta$-integral}
if $A \subset \Delta$, $A \cap \Delta_0 \ne \emptyset$,
and $T A \subset \Delta^{\perp}$ over the smooth points of $A \cap \Delta_0$.

\begin{lemma}\label{lemma-integral-is-unique}
  Let $A_1$ and $A_2$ be two irreducible $\Delta$-integral subvarieties.
  Assume $\dim A_1 =\dim A_2 = \codim_X \Delta$.
  Then either $A_1=A_2$ or $A_1 \cap A_2 \subset \Delta \setminus \Delta_0$.
\end{lemma}
\noprf

\begin{thm}
  Consider a general point $x \in X$.
  Then:
  \begin{enumerate}
     \item \label{item-Dx-cointegrable}
        $D_x$ (as reduced, but possibly not irreducible subvariety of $X$) is $F$-cointe\-grable.
     \item \label{item-lines-are-Dx-integrable}
        For general $y$ in any of 
             the irreducible components of $C_x$ 
             all lines through $y$ 
             are either $D_x$-integral
             or contained in the singular locus of~$D_x$.
     \item \label{item-unique-presentation-as-sum-of-2-lines}
        For general $z$ in any of the irreducible components of $D_x$ 
        the intersection $C_x \cap C_z$ is a unique point 
        and the chain of two lines connecting $x$ to $z$ is unique.
  \end{enumerate}
\end{thm}

\begin{prf}
   Part \ref{item-Dx-cointegrable} is immediate, since $D_x$ is a divisor, 
   by Theorem~\ref{thm-Dx-is-a-divisor}.
   
   To prove part \ref{item-lines-are-Dx-integrable} let $l$ be a line through $y$ which is not contained in the singular locus of $D_x$.
   Then by Proposition~\ref{prop-tangent-to-Dx}:
    \[
      T_z D_x \cap F_{z} = (F|_{l})^{\ge 0}
    \]
    and  for general $z \in l$ we have $(T_z D_x \cap F_{z})^{\perp_{\ud \theta_z}} \subset F_{z}$ is the $\ccO(2)$ part, 
    i.e.~the part tangent to $l$. So $l$ is $D_x$-integral as claimed.

    To prove \ref{item-unique-presentation-as-sum-of-2-lines}, 
    let $U \subset X$ be an open dense subset of points $u \in X$ where two different lines through $u$ 
    do not share the tangent direction and do not meet in any other point.
    Note that since $x$ is a general point, $x \in U$ 
    and thus each irreducible component of $C_x$ and $D_x$ intersects $U$.
    Thus generality of $z$  implies that $z \in U$
    and thus each irreducible component of $C_z$ and $D_z$ intersects $U$.
    Also $C_x \cap C_z$ intersects $U$. So fix $y \in C_x \cap C_z \cap U$.

    By \ref{item-lines-are-Dx-integrable} and Lemma \ref{lemma-integral-is-unique}
    the line $l_z$ through $z$ which intersects $C_x$ is unique.
    In the same way let $l_x$ be the unique line through $x$ intersecting $C_z$.
    Thus 
    \[
      C_x \cap C_z = l_x \cap l_z.
    \]
    In particular, $y \in  l_x \cap l_z$.
    But since $y \in U$ the intersection $l_x \cap l_z$ is just one point and therefore:
    \[
      C_x \cap C_z = \set{y}.
    \]
\end{prf}

As a consequence of part \ref{item-unique-presentation-as-sum-of-2-lines} of the theorem 
the surjective map $\pi_{13}\colon C^3 \to D$ is birational.
Thus consider the inverse rational map $D \dashrightarrow C^3$ and compose it with
the projection on the middle coordinate $\pi_2: C^3 \to X$.
We define the composition to be the \emph{bracket map}:
\[
  [\cdot, \cdot]^D \colon D \dashrightarrow C^3 \stackrel{\pi_2}{\to} X.
\]
In this setting, for $(x,z) \in D$, one has $[x,z]^D = y = C_x \cap C_z$,
 whenever the intersection is just one point.

\begin{thm}
    If $X$ is the adjoint variety of $G$, then the bracket map defined above agrees with the Lie bracket on 
    $\gotg$, in the following sense:
    Let  $\xi, \zeta \in \gotg$ and set $\eta := [\xi,\zeta]$ (the Lie bracket on $\gotg$).
    Denote by $x$, $y$ and  $z$ the projective classes in $\PP \gotg$
    of $\xi$, $\eta$ and $\zeta$ respectively.
    If $x \in D_z$ and $\eta \ne 0$,
    then the bracket map satisfies $[x,z]^D=y$. 
\end{thm}

\begin{prf}
    It is enough to prove the statement for a general pair $(x,z) \in D$.
    Suppose further $w \in C_z$ is a general point.
    Then the pair $(x,w) \in X \times X$ is a general pair.
    Thus by Theorem~\ref{thm-grading-for-homogeneous},
    we may assume $\xi \in \gotg_2$ and $\zeta \in \gotg_{-1}$.
    The restriction of the Lie bracket to $[\xi , \gotg_{-1}]$
    determines an isomorphism  $\gotg_{-1} \to \gotg_{1}$ 
    of $\gotg_0$-modules.
    In particular the minimal orbit $X \cap \PP\gotg_{-1}$ 
    is mapped onto $X \cap \PP\gotg_{1}$ under this isomorphism.
    In particular $y \in X \cap \PP\gotg_{1} \subset C_x$
    (see Lemma~\ref{lemma-homogeneous-then-X-cap-g1-subset-Cx}).
    Analogously $y \in C_z$, so $y \in C_x \cap C_z$.
\end{prf}

% \bibliography{references}
\bibliography{contact_ujednolicone}

\bibliographystyle{alpha_four}%%--%%  \end{document}

\end{document}